\documentclass[reqno,12pt]{amsart}

\theoremstyle{remark}

\numberwithin{equation}{section}

\newcommand{\al}{\alpha}
\newcommand{\be}{\beta}
\newcommand{\lm}{\lambda}
\newcommand{\s}{\sigma}
\newcommand{\ta}{\theta}

\newcommand{\half}{\frac 12}
\newcommand{\ty}{\infty}

\allowdisplaybreaks

\begin{document}

\title[Formulas for multibasic theta hypergeometric series]
{Summation, transformation,\\
and expansion formulas for\\
multibasic theta hypergeometric series}

\author{George Gasper}
\address{Department of Mathematics, Northwestern University,
2033 Sheridan Road, Evanston, IL 60208-2730, USA}
\email{george@math.northwestern.edu}
\urladdr{http://www.math.northwestern.edu/{\textasciitilde}george}

\author[Michael Schlosser]{Michael Schlosser$^*$}
\address{Institut f\"ur Mathematik der Universit\"at Wien,
Nordbergstra{\ss}e 15, A-1090 Wien, Austria}
\email{schlosse@ap.univie.ac.at}
\urladdr{http://www.mat.univie.ac.at/{\textasciitilde}schlosse}

\thanks{$^*$The seconds author's research was supported by
an APART fellowship of the Austrian Academy of Sciences,
by FWF Austrian Science Fund grant \hbox{P17563-N13},
and by EC's IHRP Programme, grant HPRN-CT-2001-00272
``Algebraic Combinatorics in Europe''.}

\begin{abstract}
After reviewing some fundamental facts from the theory of theta
hypergeometric series we derive, using indefinite summation,
several summation, transformation, and expansion formulas for
multibasic theta hypergeometric series. Some of the identities
presented here generalize corresponding formulas given in Chapter~11
of the Gasper and Rahman book [{\em Basic hypergeometric series},
$2^{nd}$ ed., Encyclopedia of Mathematics And Its Applications~96,
Cambridge University Press, Cambridge, 2004].
\end{abstract}

\maketitle

\section{Introduction}

By convention, a series $\sum u_n$ is called a hypergeometric
series if $g(n) = u_{n+1}/u_n$ is a rational function of $n$. It is
called a $q$- (or basic) hypergeometric series if $g(n)$ is a rational
function of $q^n$. More generally, such a series is called an 
{\it elliptic hypergeometric series} if $g(n)$ is an elliptic
(doubly periodic meromorphic) function of $n$ with $n$ considered
as a complex variable. We refer the reader to \cite[Ch.~11]{GR}
for some motivation for considering these three classes of series
and to  \cite{GR} in general for a treatise
on basic hypergeometric series.

In a path-breaking paper, Frenkel and Turaev~\cite{FT2}
in their work on elliptic $6j$-symbols (introduced by
Date {\it et al.}~\cite{D1} as elliptic solutions of the
Yang-Baxter equation \cite{B1,B2}) introduced elliptic analogues
of very-well-poised basic hypergeometric series.
In particular, using the tetrahedral symmetry of the elliptic
$6j$-symbols and the finite dimensionality of cusp forms, 
they derived elliptic analogues of  Bailey's transformation
formula (cf.\ \cite[Eq.~(2.9.1)]{GR}) for terminating
$_{10}\phi_9$ series and of Jackson's $_8\phi_7$
summation formula (cf.\ \cite[Eq.~(2.6.2)]{GR}).
Elliptic hypergeometric series and their extensions to
theta hypergeometric series became an increasingly active
area of research (see \cite[Sec.~11.1]{GR} for some references).
So far, many formulas for very-well-poised basic hypergeometric
series have already been extended to the elliptic setting.
Some formulas for {\em multibasic} elliptic hypergeometric series
appeared in work of Warnaar~\cite{W}. Here we consider yet other
identities involving multiple bases and theta functions,
special cases of which have already been presented in \cite{GR}.

We start in Section~\ref{secellser} with the elliptic shifted
factorials, Spiridonov's~\cite{S1} $_{r+1}E_r$ theta hypergeometric
series notation and its very-well-poised $_{r+1}V_r$ special case,
and then point out some of their main properties.  
We also present the Frenkel and Turaev summation and
transformation formulas. In Section~\ref{secmult} we derive
theta hypergeometric extensions of some of the summation
and transformation formulas in \cite[Secs.~3.6--3.8]{GR}.
To give just one example, here is a transformation formula for a
``split-poised'' theta hypergeometric $_{12}E_{11}$ series
\begin{align}\label{splitpoised}
&_{12}E_{11} \big(a, qa^{\half}, -qa^{\half}, q a^{\half}/p^{\half},
-q a^{\half}p^{\half}, b, c, a/bc, q^{-n}, B/Aq^n, C/Aq^n, 1/BCq^n;\cr
& \qquad\quad a^{\half}, -a^{\half}, a^{\half}p^{\half}, -a^{\half}/p^{\half},
aq/b, aq/c, bcq, 1/Aq^n, 1/Bq^n, 1/Cq^n, BC/Aq^n;\cr
&\hspace{30.4em} q,p;-1\big)\cr
&=\frac{(aq, bq, cq, aq/bc, Aq/B, Aq/C, BCq; q,p)_n}
{(Aq, Bq, Cq, Aq/BC, aq/b, aq/c, bcq;q,p)_n}\cr
&\times{}_{12}E_{11} \big(A, qA^{\half}, -qA^{\half}, q A^{\half}/p^{\half},
-q A^{\half}p^{\half}, B, C, A/BC, q^{-n}, b/aq^n, c/aq^n, 1/bcq^n;\cr
& \qquad\quad\;\;A^{\half}, -A^{\half}, A^{\half}p^{\half},
-A^{\half}/p^{\half}, Aq/B, Aq/C, BCq, 1/aq^n, 1/bq^n, 1/cq^n, bc/aq^n;\cr
&\hspace{31.4em} q,p;-1\big)
\end{align}
for $n=0,1,\ldots$ (see \eqref{spser} for the notation),
a recast of \eqref{splitpoised1}, which extends the
transformation formula for a split-poised $_{10}\phi_9$ series
given in \cite[Ex.$\,$3.21]{GR}.
Most of these extensions have recently been presented in 
\cite[Ch.~11]{GR},
where reference was made to an earlier (2003) version of this paper.
However, some formulas in \cite{GR} (in particular, in Sec. 11.6
and Exercises 11.25--11.26) have been further generalized in the current
version of this paper. 
The selection of formulas we give is by no means
exhaustive, but they do serve to illustrate some of the possibilities
for deriving summations, transformations and expansions for multibasic 
theta hypergeometric series. We wish to thank Mizan Rahman and
Ole Warnaar for some helpful correspondences.

\section{Elliptic and theta hypergeometric series}\label{secellser}

As in \cite{GR} we define a modified Jacobi theta function with argument
$x$ and nome $p$ by
\begin{equation}\label{mjtf}
\ta (x; p) = (x, p/x; p)_\ty =(x; p)_\ty (p/x; p)_\ty\,,\quad\quad
\ta (x_1, \ldots, x_m;
p) = \prod^m_{k=1} \ta (x_k; p),
\end{equation}
where $ x, x_1, \ldots, x_m \ne 0,\ |p| < 1,$ and $(x; p)_\ty=
\prod^\infty_{k=0}(1-x p^k) .$
Also, following Warnaar~\cite{W}, we define an {\it elliptic} (or {\it theta})
{\it shifted factorial} analogue of the $q$-shifted factorial  by
\begin{equation}\label{defepoch}
(a;q,p)_n = \begin{cases}
\prod^{n-1}_{k=0} \ta (aq^k;p),& n = 1, 2, \ldots\,,\cr
1,& $n = 0,$\cr
1/\prod^{-n-1}_{k=0} \ta (aq^{n+k};p),& n = -1, -2, \ldots,
\end{cases}
\end{equation}
and let
\begin{equation*}
(a_1, a_2, \ldots, a_m;q, p)_n = \prod^m_{k=1} (a_k;q,p)_n,
\end{equation*}
where $a, a_1,\ldots,a_m \neq 0$.
Notice that $\ta (x;0) = 1-x$ and, hence, $(a;q, 0)_n = (a;q)_n$ is a
$q$-{\it shifted factorial} in base $q.$ Thus, the parameters
 $q$ and $p$ in $(a;q,p)_n$ are called the {\it base} and {\it nome},
respectively, and
 $(a;q,p)_n$ is called the $q,p$-{\it shifted factorial}.
A list of useful identities for manipulating the
$q,p$-shifted factorials (and related objects such as
$q,p$-{\it binomial coefficients}, or
{\it elliptic binomial coefficients}, and {\it elliptic gamma functions})
is given in  \cite[Sec.~11.2]{GR}.

Following Spiridonov~\cite{S1},
 an $_{r+1}E_r$ {\it theta hypergeometric series}
with base $q$ and nome $p$ is formally defined by
\begin{align}\label{spser}
&_{r+1}E_r (a_1, a_2, \ldots, a_{r+1}; b_1, \ldots, b_r;q, p;z)\cr
&=\sum^\ty_{n=0} \frac{(a_1, a_2, \ldots, a_{r+1};q, p)_n}
{(q, b_1, \ldots, b_r;q,p)_n}z^n,
\end{align}
where, as usual, it is assumed that the parameters are such
that each term in the  series is well-defined. 
If $z$ and the $a$'s and $b$'s are independent of
$p$, then it follows that
\begin{align*}
&\lim_{p\to 0} {}_{r+1}E_r(a_1, \ldots, a_{r+1}; b_1, \ldots, b_r; q,p;z)\cr
&=\, {}_{r+1}E_r(a_1, \ldots, a_{r+1}; b_1, \ldots, b_r; q,0;z)\cr
&=\, {}_{r+1}\phi_r (a_1, \ldots, a_{r+1}; b_1, \ldots, b_r; q,z),
\end{align*}
where the limit of the  series is a termwise limit.
See \cite[Sec.~11.2]{GR} for more details and a discussion of convergence
of the series in \eqref{spser}.

As in \cite{S1},  a (unilateral or bilateral) series
$\sum c_n$ is called an {\it elliptic hypergeometric series} if
$g(n) = c_{n+1}/c_n$ is an elliptic function of $n$ with $n$
considered as a complex variable; i.e., the function $g(x)$
is a doubly periodic
meromorphic function of the complex variable $x$. For the 
$_{r+1}E_r$ series in \eqref{spser} it is clear that
\begin{equation*}
g(x) = z\prod^{r+1}_{k=1}\frac{\ta (a_kq^x;p)}
{\ta (b_k q^x;p)}
\end{equation*}
with $b_{r+1} = q$.
It is not difficult to show (see \cite{GR}) that when
\begin{equation}\label{ebalc}
a_1a_2\dots a_{r+1}=(b_1b_2\dots b_r)q,
\end{equation}
$g(x)$ is an elliptic
(i.e., doubly periodic meromorphic) function of $x$.
Therefore, \eqref{ebalc} is called the
{\it elliptic balancing condition}, and  $_{r+1}E_r$ is said to be
{\it elliptically balanced} ({\it E-balanced}) when \eqref{ebalc}
holds.

Corresponding to the basic hypergeometric special case (cf.\ \cite{GR}),
 the $_{r+1}E_r$ series in \eqref{spser} is called {\it well-poised} if
\begin{equation}\label{qab}
qa_1 = a_2 b_1 = a_3 b_2 = \ldots = a_{r+1} b_r,
\end{equation}
in which case we find that the elliptic balancing condition \eqref{ebalc}
reduces to
\begin{equation*}
a^2_1a^2_2\cdots a^2_{r+1} = (a_1 q)^{r+1}.
\end{equation*}
Using \eqref{defepoch} we see that
\begin{equation*}
\frac{\ta (a q^{2n};p)}
{\ta (a;p)} = \frac{(qa^{\half}, -qa^{\half},
        q a^{\half}/p^{\half}, -q a^{\half}p^{\half};q,p)_n}
{(a^{\half}, -a^{\half}, a^{\half}p^{\half},
        - a^{\half}/p^{\half};q,p)_n} (-q)^{-n}
\end{equation*}
is an elliptic analogue of the quotient
\begin{equation*}
\frac{1-aq^{2n}}{1-a} = \frac{(qa^{\half}, -qa^{\half};q)_n}
{(a^{\half}, - a^{\half};q)_n},
\end{equation*}
which is clearly the very-well-poised part of the very-well-poised
basic hypergeometric $_{r+1}W_r$ series in \cite[Eq.~(2.1.11)]{GR}.
Hence, the $_{r+1}E_r$ series in \eqref{spser} is called
{\it very-well-poised} when it is well-poised, $r\ge 4$, and
\begin{equation}\label{q2r1}
a_2 = qa^{\half}_1, \ a_3 = -qa^{\half}_1,\ a_4 =
q a_1^{\half}/p^{\half}, \ a_5 = -q a_1^{\half}p^{\half}.
\end{equation}

Analogous to Spiridonov~\cite[Eq.~(2.15)]{S1}, an
$_{r+1}V_r$ very-well-poised theta hypergeometric series is defined by
\begin{align}\label{vser}
&{}_{r+1}V_r (a_1; a_6, a_7, \ldots, a_{r+1}; q,p;z)\cr
&=\sum^\ty_{n=0}\frac{\ta (a_1 q^{2n};p)}{\ta (a_1; p)}
\frac{(a_1, a_6, a_7,
\ldots, a_{r+1};q, p)_n}
{(q, a_1 q/a_6, a_1 q/a_7, \ldots, a_1q/a_{r+1};q,p)_n}\,(q z)^n.
\end{align}
Thus, if \eqref{qab} and \eqref{q2r1} hold, then
\begin{align*}
&{}_{r+1}V_r (a_1; a_6, a_7, \ldots, a_{r+1}; q,p;z)\cr
&= {}_{r+1}E_r (a_1, a_2, \ldots, a_{r+1}; b_1, \ldots, b_r;q,
p;-z),
\end{align*}
and the ${}_{r+1}V_r$ series is elliptically
balanced if and only if
\begin{equation*}
(a^2_6a^2_7\cdots a^2_{r+1}) q^2= (a_1 q)^{r-5}.
\end{equation*}
If the argument $z$ in the ${}_{r+1}V_r$ series equals $1,$ then
we suppress it and denote the series in \eqref{vser} by the simpler
notation ${}_{r+1}V_r (a_1; a_6, a_7, \ldots, a_{r+1}; q,p).$
When the parameters $a_1, a_6, a_7, \ldots, a_{r+1}$
are independent of $p,$
\begin{align*}
&\lim_{p\to 0} {}_{r+1}V_r (a_1; a_6, a_7, \ldots, a_{r+1}; q,p)\cr
&=\, {}_{r-1}W_{r-2}(a_1; a_6, \ldots, a_{r+1}; q,q),
\end{align*}
from which it follows  that there is a shift $r\to r-2$ when taking the
$p \to 0$ limit, and that the $p\to 0$ limit of a
${}_{r+1}V_{r}(a_1;a_6,a_7,\ldots,a_{r+1};q,p)$ series with
$a_1, a_6, a_7, \ldots, a_{r+1}$ independent of $p$ is a
${}_{r-1}W_{r-2}$ series.

Frenkel and Turaev~\cite{FT2} showed the following elliptic analogue of
Bailey's $_{10}\phi_{9}$ transformation formula~\cite[Eq.~(2.9.1)]{GR}
\begin{align}\label{ft1211}
&_{12}V_{11} (a; b, c, d, e, f, \lm aq^{n+1}/ef, q^{-n};q,p)\cr
&=\frac{(aq, aq/ef, \lm q/e, \lm q/f;q,p)_n}
{(aq/e, aq/f, \lm q/ef, \lm q;q,p)_n}\cr
&\quad \times {}_{12}V_{11} (\lm; \lm b/a, \lm c/a, \lm d/a, e,
f, \lm aq^{n+1}/ef, q^{-n};q,p)
\end{align}
for $n=0,1,\ldots,$ provided that the  balancing condition
\begin{equation}\label{econd}
bcdef(\lm aq^{n+1}/ef)q^{-n}q= (a q)^3,
\end{equation}
which is clearly equivalent to $\lm= qa^2/bcd,$ holds.
Notice that each of the series in \eqref{ft1211}
is E-balanced when \eqref{econd} holds.
If we set $\lm =a/d$ in \eqref{ft1211}, then we obtain a summation formula
for a $_{10}V_9$ series which is an elliptic analogue of Jackson's
$_8\phi_7$ summation formula~\cite[Eq.~(2.6.2)]{GR}
and of Dougall's $_7F_6$ summation formula~\cite[Eq.~(2.1.6)]{GR}.
After a change in parameters,
this summation formula can be written in the form:
\begin{equation}\label{ft109}
{}_{10}V_9 (a; b, c, d, e, q^{-n};q,p) =
\frac{(aq, aq/bc, aq/bd, aq/cd;q,p)_n}
{(aq/b, aq/c, aq/d, aq/bcd;q,p)_n}
\end{equation}
for $n=0,1,\ldots,$ provided that the elliptic balancing condition
$bcde=a^2q^{n+1},$ which can be written in the form
\begin{equation}
(bcdeq^{-n})q =(aq)^2,
\end{equation}
holds. It is obvious that if $a, b, c, d, e$
are independent of $p$, then \eqref{ft109} tends to Jackson's 
$_8\phi_7$ summation formula~\cite[Eq.~(2.6.2)]{GR} as $p\to 0$.
For a further discussion of \eqref{ft1211} and \eqref{ft109}
including different proofs, see
Gasper and Rahman~\cite[Secs.~11.2, 11.4, 11.5]{GR}.

\section{Multibasic summation and transformation formulas
for theta hypergeometric series}\label{secmult}

We first observe that if the parameter $a$ in \eqref{ft109} 
is replaced by $a/q,$ then it follows
that the $n=1$ case of \eqref{ft109} is equivalent to the identity
\begin{equation}\label{ft109n1}
1 - \frac{\ta (b, c, d, a^2/bcd;p)}
{\ta (a/b, a/c, a/d, bcd/a; p)} = \frac{\ta (a, a/bc, a/bd, a/cd;p)}
{\ta (a/bcd, a/d, a/c, a/b;p)}.
\end{equation}
More generally, by replacing $a$ in \eqref{ft1211} by $a/q$ it follows
that the $n=1$ case of \eqref{ft1211} is equivalent to the identity
\begin{align}\label{ft1211n1}
&1 - \frac{\ta (b, c, d, e, f,g;p)}
{\ta (a/b, a/c, a/d, a/e, a/f, a/g; p)} \cr
&=\frac{\ta (a, a/ef, a^2/bcde, a^2/bcdf;p)}
{\ta (a^2/bcdef, a^2/bcd, a/f, a/e;p)} \cr
&\quad \times\left[1-\frac{\ta (a/bc, a/bd, a/cd, e, f, g;p)}
{\ta (a/d, a/c, a/b, a^2/bcde, a^2/bcdf, a^2/bcdg;p)}\right]
\end{align}
with $a^3=bcdefg$, which is equivalent to the identity in
\cite[Ex.~5.22]{GR}. Next, define
\begin{equation}\label{gprod}
\prod^n_{k=m} a_k = \begin{cases}
a_m a_{m+1} \cdots a_n,& m \le n,\cr
1,& m = n+1,\cr
(a_{n+1} a_{n+2} \cdots a_{m-1})^{-1},& m \ge n + 2,
\end{cases}
\end{equation}
for $n,m = 0$, $\pm 1$, $\pm 2, \ldots\,$, and let
\begin{equation}\label{Un}
U_n = \prod^{n-1}_{k=0}\frac{\ta (b_k, c_k, d_k, e_k, f_k, g_k;p)}
{\ta (a_k/b_k, a_k/c_k, a_k/d_k, a_k/e_k, a_k/f_k, a_k/g_k;p)}
\end{equation}
where $a^3_k = b_k c_k d_k e_k f_k g_k$ for
$k = 0$, $\pm 1, \pm 2, \ldots\,$, and it is assumed that the
$a$'s, $b$'s, $c$'s, $d$'s, $e$'s, $f$'s, $g$'s, are complex numbers
such that $U_n$ is well defined for $n = 0, \pm 1, \pm 2, \ldots$ . 
Now use \eqref{ft1211n1} with $a$, $b$, $c$, $d$, $e$, $f$, $g$
replaced by $a_k$, $b_k$, $c_k$, $d_k$, $e_k$, $f_k$, $g_k$ respectively,
to get the indefinite summation formula
\begin{align}\label{indefsum}
&U_{-m} - U_{n+1} = \sum^n_{k=-m} (U_k - U_{k+1})\cr
&= \sum^n_{k=-m}\frac{\ta (a^2_k/b_k c_k d_k f_k,
a^2_k/b_k c_k d_k e_k, a_k, a_k/e_k f_k;p)}
{\ta (a_k/e_k, a_k/f_k,
a^2_k/b_k c_k d_k e_k f_k, a^2_k/b_k c_k d_k;p)} U_k\cr
&\times \bigg[ 1 -\frac{\ta (a_k/c_k d_k, a_k/b_kd_k,
a_k/b_k c_k, e_k, f_k, g_k;p)}
{\ta (a_k/b_k, a_k/c_k, a_k/d_k, a^2_k/b_k c_k d_k e_k,
 a^2_k/b_k c_k d_k f_k, a^2_k/b_k c_k d_k g_k;p)}\bigg]
\end{align}
for $n,m = 0$, $\pm 1, \pm 2, \ldots$,  where
$a^3_k = b_k c_k d_k e_k f_k g_k$ for $k = 0$, $\pm 1, \pm 2, \ldots\,$.
Since $U_0 = 1$ by \eqref{gprod} and $\ta (a_k/c_kd_k;p)=0$ when
$a_k=c_kd_k$, setting $m = 0$ and $a_k=c_kd_k$ for $k=0,1,\dots,n$
in \eqref{indefsum} yields after relabelling the summation formula
\begin{align}\label{sumf}
&\sum^n_{k=0} \frac{\ta (a_k, a_k/b_k c_k, a_k/b_k d_k, a_k/c_k d_k;p)}
{\ta (a_k/b_k c_k d_k, a_k/d_k, a_k/c_k, a_k/b_k;p)}\cr
&\quad\times \prod^{k-1}_{j=0}\frac{\ta (b_j, c_j, d_j, a^2_j/b_j c_j d_j;p)}
{\ta (a_j/b_j, a_j/c_j, a_j/d_j, b_j c_jd_j/a_j;p)}\cr
&= 1 - \prod^n_{j=0}\frac{\ta (b_j, c_j, d_j, a^2_j/b_jc_j d_j;p)}
{\ta (a_j/b_j, a_j/c_j, a_j/d_j, b_j c_j d_j/a_j;p)}
\end{align}
for $n = 0, 1, \ldots\,$,
which is equivalent to Warnaar's formula \cite[Eq.~(3.2)]{W}. 
When $p = 0$ the above formula reduces to a summation formula of Macdonald
that was first published in Bhatnagar and Milne~\cite[Thm.~2.27]{BM},
and  contains the summation formulas by W.~Chu~\cite[Thms.~A, B, C]{C}
as special cases.

Observe that in \eqref{ft109n1}, \eqref{ft1211n1}, \eqref{Un},
\eqref{indefsum} and \eqref{sumf} the
components of each quotient of products of theta functions 
 have been arranged so that the 
well-poised property of these quotients is clearly displayed;
e.g., in the second sum in \eqref{indefsum} the quotient of the theta
functions in front of $U_k$ is arranged so that each product of
corresponding numerator and denominator parameters equals
$a^3_k/b_kc_kd_ke_kf_k$, and each of the corresponding products
in the quotient of theta functions inside the square bracket equals
$a^2_k/b_kc_kd_k$.

If we let
\begin{equation*}
a_k = aw^k, \ b_k = bq^k, \ c_k = cr^k, \ d_k = ds^k,
\ e_k = et^k, \ f_k = fu^k, \ g_k = gv^k,
\end{equation*}
with $a^3=bcdefg$ and $w^3=qrstuv$, then the product $U_n$ reduces to
\begin{align*}
\tilde U_n =&\frac{(b;q,p)_n (c;r,p)_n (d;s,p)_n}
{(a/b;w/q,p)_n (a/c;w/r,p)_n (a/d;w/s,p)_n}\cr
&\times\frac{(e;t,p)_n (f;u,p)_n (g;v,p)_n}
{(a/e;w/t,p)_n (a/f;w/u,p)_n (a/g;w/v,p)_n}
\end{align*}
and, by applying \eqref{indefsum} and  some
elementary identities for $q,p$-shifted factorials
(listed in \cite[Sec.~11.2]{GR}), we obtain  the following
indefinite multibasic theta hypergeometric summation formula
\begin{align}\label{indm}
&\sum^n_{k=-m}
\frac{\ta (aw^k, a(w/tu)^k/ef, fg(uv/w)^k/a, eg(tv/w)^k/a;p)}
{\ta (g(v/w)^k/a, efg(tuv/w)^k/a, a(w/u)^k/f, a(w/t)^k/e;p)}\cr
&\quad\times\frac{(b;q,p)_k (c;r,p)_k (d;s,p)_k}
{(a/b;w/q,p)_k (a/c;w/r,p)_k (a/d;w/s,p)_k}\cr
&\quad\times\frac{(e;t,p)_k (f;u,p)_k (g;v,p)_k}
{(a/e;w/t,p)_k (a/f;w/u,p)_k (a/g;w/v,p)_k}\cr
&\quad\times\bigg[1-\frac{\ta (a(w/rs)^k/cd, a(w/qs)^k/bd, a(w/qr)^k/bc;p)}
{\ta (a(w/q)^k/b, a(w/r)^k/c, a(w/s)^k/d;p)}\cr
&\quad\qquad\times
\frac{\ta (et^k, fu^k, gv^k;p)}
{\ta (fg(uv/w)^k/a, eg(tv/w)^k/a, ef(tu/w)^k/a;p)}\bigg]\cr
&=\frac{(bw/aq;w/q,p)_m (cw/ar;w/r,p)_m (dw/as;w/s,p)_m}
{(q/b;q,p)_m (r/c;r,p)_m (s/d;s,p)_m}\cr
&\quad\times\frac{(ew/at;w/t,p)_m (fw/au;w/u,p)_m (gw/av;w/v,p)_m}
{(t/e;t,p)_m (u/f;u,p)_m (v/g;v,p)_m}\cr
&\quad-\frac{(b;q,p)_{n+1} (c;r,p)_{n+1} (d;s,p)_{n+1}}
{(a/b;w/q,p)_{n+1} (a/c;w/r,p)_{n+1} (a/d;w/s,p)_{n+1}}\cr
&\qquad \times\frac{(e;t,p)_{n+1} (f;u,p)_{n+1} (g;v,p)_{n+1}}
{(a/e;w/t,p)_{n+1} (a/f;w/u,p)_{n+1} (a/g;w/v,p)_{n+1}}
\end{align}
for $n,m = 0, \pm 1, \pm 2, \ldots,$ where $a^3=bcdefg$ and $w^3=qrstuv$.

If we set $p = 0$ and assume that
\begin{equation*}
\max (|q|, |r|, |s|, |t|, |u|, |v|,
|w/q|, |w/r|, |w/s|, |w/t|, |w/u|, |w/v|) < 1,
\end{equation*}
then letting $n$ or $m$ in \eqref{indm}
tend to infinity shows that this special 
case of \eqref{indm} also holds with $n$ and/or $m$ replaced by $\infty$,
just as in the special case \cite[Eq.~(3.6.14)]{GR}. Thus we have extended
\cite[Eq.~(3.6.14)]{GR} to the
bilateral multibasic summation formula
\begin{align}\label{indmrat}
&\sum^\ty_{k=-\ty}
\frac{(1-aw^k)(1-a(w/tu)^k/ef)(1-fg(uv/w)^k/a)(1-eg(tv/w)^k/a)}
{(1-g(v/w)^k/a)(1-efg(tuv/w)^k/a)(1-a(w/u)^k/f)(1-a(w/t)^k/e)}\cr
&\quad\times\frac{(b;q)_k (c;r)_k (d;s)_k (e;t)_k (f;u)_k (g;v)_k}
{(a/b;w/q)_k (a/c;w/r)_k (a/d;w/s)_k (a/e;w/t)_k (a/f;w/u)_k (a/g;w/v)_k}\cr
&\quad\times\bigg[1-\frac{(1-a(w/rs)^k/cd)(1-a(w/qs)^k/bd)(1-a(w/qr)^k/bc)}
{(1-a(w/q)^k/b)(1-a(w/r)^k/c)(1-a(w/s)^k/d)}\cr
&\quad\qquad\times
\frac{(1-et^k)(1-fu^k)(1-gv^k)}
{(1-fg(uv/w)^k/a)(1-eg(tv/w)^k/a)(1-ef(tu/w)^k/a)}\bigg]\cr
&=\frac{(bw/aq;w/q)_\ty (cw/ar;w/r)_\ty (dw/as;w/s)_\ty}
{(q/b;q)_\ty (r/c;r)_\ty (s/d;s)_\ty}\cr
&\quad\times\frac{(ew/at;w/t)_\ty (fw/au;w/u)_\ty (gw/av;w/v)_\ty}
{(t/e;t)_\ty (u/f;u)_\ty (v/g;v)_\ty}\cr
&\quad-\frac{(b;q)_\ty (c;r)_\ty (d;s)_\ty (e;t)_\ty (f;u)_\ty (g;v)_\ty}
{(a/b;w/q)_\ty (a/c;w/r)_\ty (a/d;w/s)_\ty
(a/e;w/t)_\ty (a/f;w/u)_\ty (a/g;w/v)_\ty},
\end{align}
where $a^3=bcdefg$ and $w^3=qrstuv$, and
\begin{equation*}
\max (|q|, |r|, |s|, |t|, |u|, |v|,
|w/q|, |w/r|, |w/s|, |w/t|, |w/u|, |w/v|) < 1.
\end{equation*}

Even though we cannot let $n \to \ty$ or $m \to \ty$ in \eqref{indm}
when $p\ne 0$ to derive summation formulas 
for nonterminating theta hypergeometric series (because 
$\lim_{a\to 0} \ta (a;p)$ does not exist when $p\ne 0$), it is possible 
in some special cases to let 
$n\to \ty$ or $m\to\ty$ in \eqref{indefsum} to obtain summation
formulas for nonterminating series containing products of 
certain theta functions. In particular, if we
denote the $k$th factor in the product representation \eqref{Un} for $U_n$ by
\begin{equation*}
z_k = \frac{\ta (b_k, c_k, d_k, e_k, f_k, a_k^3/b_k c_k d_k e_k f_k; p)}
{\ta (a_k/b_k, a_k/c_k, a_k/d_k, a_k/e_k, a_k/f_k, b_kc_kd_ke_kf_k/a_k^2;p)}
\end{equation*}
and observe that 
\begin{equation*}
\lim_{b\to a^{{\half}}}\,\frac{\ta (b;p)}{\ta (a/b;p)} = 1,\quad |p| < 1,
\end{equation*}
when $a$ is not an integer power of $p$,
then it follows that there exist bilateral sequences of the 
$a$'s, $b$'s, $c$'s, $d$'s, $e$'s, and $f$'s in \eqref{indefsum} such that
$\Re z_k > 0$  for integer $k$ and the series
\begin{equation}\label{logcond}
\sum_{k =-\ty}^\ty \ \log z_k \ \  \ \ {\rm converges},
\end{equation}
where $\log z_k$ is the principal branch of the logarithm (choose, e.g., 
$b_k$, $c_k$, $d_k$, $e_k$, and $f_k$ so close to
$a^{{\half}}_k$ that $|\log z_k| < 1/k^2$ 
for $k = \pm 1, \pm 2, \ldots\,$).
Then both of the limits $\lim_{n\to \ty} U_n$
and $\lim_{m\to \ty} U_{-m}$ exist, and we obtain the 
bilateral summation formula (which extends \cite[Eq.~(11.6.8)]{GR})
\begin{align}
&\sum^\ty_{k=-\ty}\frac{\ta (a_k, a_k/e_k f_k, a_k^2/b_kc_kd_ke_k,
a_k^2/b_kc_kd_kf_k;p)}
{\ta (a_k^2/b_k c_k d_k e_k f_k, a_k^2/b_k c_k d_k, a_k/f_k, a_k/e_k;p)}\cr
&\quad \times \prod^{k-1}_{j=0}
\frac{\ta (b_j, c_j, d_j,  e_j, f_j, a_j^3/b_j c_j d_j e_j f_j;p)}
{\ta (a_j/b_j, a_j/c_j, a_j/d_j,
a_j/e_j, a_j/f_j, b_jc_jd_je_jf_j/a_j^2;p)}\cr
&\quad \times \bigg[1-\frac{\ta (a_k/b_kc_k, a_k/b_k d_k, a_k/c_kd_k;p)}
{\ta (a_k/d_k, a_k/c_k, a_k/b_k;p)}\cr
&\quad\qquad \times\frac{\ta(e_k, f_k, g_k;p)}
{\ta (a_k^2/b_k c_k d_k e_k, a_k^2/b_k c_k d_k f_k,
a_k^2/b_k c_k d_k g_k;p)}\bigg]\cr
&=\prod^{-1}_{k=-\ty}
\frac{\ta (a_k/b_k, a_k/c_k, a_k/d_k,
a_k/e_k, a_k/f_k, b_kc_kd_ke_kf_k/a_k^2;p)}
{\ta (b_k, c_k, d_k, e_k, f_k, a_k^3/b_k c_k d_k e_k f_k; p)}\cr
&\quad -\prod^\ty_{k=0}
\frac{\ta (b_k, c_k, d_k, e_k, f_k, a_k^3/b_k c_k d_k e_k f_k; p)}
{\ta (a_k/b_k, a_k/c_k, a_k/d_k, a_k/e_k,
a_k/f_k, b_kc_kd_ke_kf_k/a_k^2;p)}
\end{align}
with $a_k^3=b_kc_kd_ke_kf_kg_k$ for $k = 0, \pm 1, \pm 2, \ldots$, and
$a_k$, $b_k$, $c_k$, $d_k$, $e_k$, $f_k$, $g_k$ such that
\eqref{logcond} holds. 

However, it seems to be more useful to employ the patching
\begin{align*}
&\ta(a(w/t)^k/e, a(w/u)^k/f, g(v/w)^k/a;p)
(a/e;w/t,p)_k (a/f;w/u,p)_k (a/g;w/v,p)_k \cr
&=\ta(a/e, a/f, g/a;p)
(aw/et;w/t,p)_k (aw/fu;w/u,p)_k (aw/gv;w/v,p)_k (v/w)^k,
\end{align*}
to convert the $m = 0$ case of \eqref{indm} into the form
\begin{align}\label{m0}
&\sum^n_{k=0}
\frac{\ta (aw^k, a(w/tu)^k/ef, fg(uv/w)^k/a, eg(tv/w)^k/a, efg/a;p)}
{\ta (a, a/ef, fg/a, eg/a, efg(tuv/w)^k/a;p)}\cr
&\quad\times\frac{(b;q,p)_k (c;r,p)_k (d;s,p)_k}
{(a/b;w/q,p)_k (a/c;w/r,p)_k (a/d;w/s,p)_k}\cr
&\quad\times\frac{(e;t,p)_k (f;u,p)_k (g;v,p)_k}
{(aw/et;w/t,p)_k (aw/fu;w/u,p)_k (aw/gv;w/v,p)_k} (w/v)^k \cr
&\quad\times\bigg[1-\frac{\ta (a(w/rs)^k/cd, a(w/qs)^k/bd, a(w/qr)^k/bc;p)}
{\ta (a(w/q)^k/b, a(w/r)^k/c, a(w/s)^k/d;p)}\cr
&\quad\qquad\times
\frac{\ta (et^k, fu^k, gv^k;p)}
{\ta (fg(uv/w)^k/a, eg(tv/w)^k/a, ef(tu/w)^k/a;p)}\bigg]\cr
&=\frac{\ta(a/e, a/f, g/a, efg/a;p)}{\ta(eg/a, fg/a, a, a/ef;p)}\cr
&\quad\times\bigg[1-\frac{(b;q,p)_{n+1} (c;r,p)_{n+1} (d;s,p)_{n+1}}
{(a/b;w/q,p)_{n+1} (a/c;w/r,p)_{n+1} (a/d;w/s,p)_{n+1}}\cr
&\qquad \times\frac{(e;t,p)_{n+1} (f;u,p)_{n+1} (g;v,p)_{n+1}}
{(a/e;w/t,p)_{n+1} (a/f;w/u,p)_{n+1} (a/g;w/v,p)_{n+1}}\bigg]
\end{align}
where $a^3=bcdefg$ and $w^3=qrstuv$, and then to let $g=v^{-n}$ to obtain the
following multibasic theta hypergeometric generalization of
\cite[Eq.~(3.6.16)]{GR}
\begin{align}\label{csn}
&\sum^n_{k=0}\frac{\ta (aw^k, a(w/tu)^k/ef, a^2(uv/w)^k/bcde,
a^2(tv/w)^k/bcdf, a^2/bcd;p)}
{\ta (a, a/ef, a^2/bcde, a^2/bcdf, a^2(tuv/w)^k/bcd;p)}\cr
&\quad\times\frac{(b;q,p)_k (c;r,p)_k (d;s,p)_k}
{(a/b;w/q,p)_k (a/c;w/r,p)_k (a/d;w/s,p)_k}\cr
&\quad\times\frac{(e;t,p)_k (f;u,p)_k (v^{-n};v,p)_k}
{(aw/et;w/t,p)_k (aw/fu;w/u,p)_k (awv^{n-1};w/v,p)_k} (w/v)^k \cr
&\quad\times\bigg[1-\frac{\ta (a(w/rs)^k/cd, a(w/qs)^k/bd, a(w/qr)^k/bc;p)}
{\ta (a(w/q)^k/b, a(w/r)^k/c, a(w/s)^k/d;p)}\cr
&\quad\qquad\times
\frac{\ta (et^k, fu^k, v^{k-n};p)}
{\ta (f(u/w)^kv^{k-n}/a, e(t/w)^kv^{k-n}/a, ef(tu/w)^k/a;p)}\bigg]\cr
&=\frac{\ta(a/e, a/f, v^{-n}/a, efv^{-n}/a;p)}
{\ta(ev^{-n}/a, fv^{-n}/a, a, a/ef;p)},
\end{align}
where $a^3v^n=bcdef$ and $w^3=qrstuv$, and $n=0,1,\dots$ .
By letting $f\to a$ in \eqref{csn} we obtain
\begin{align}\label{ftoa}
&\sum^n_{k=0}\frac{\ta (aw^k, (w/tu)^k/e, v^{-n}(uv/w)^k,
a(tv/w)^k/bcd, a^2/bcd;p)}
{\ta (a, 1/e, a^2/bcde, a/bcd, a^2(tuv/w)^k/bcd;p)}\cr
&\quad\times\frac{(b;q,p)_k (c;r,p)_k (d;s,p)_k}
{(a/b;w/q,p)_k (a/c;w/r,p)_k (a/d;w/s,p)_k}\cr
&\quad\times\frac{(e;t,p)_k (a;u,p)_k (v^{-n};v,p)_k}
{(aw/et;w/t,p)_k (w/u;w/u,p)_k (awv^{n-1};w/v,p)_k} (w/v)^k \cr
&\quad\times\bigg[1-\frac{\ta (a(w/rs)^k/cd, a(w/qs)^k/bd, a(w/qr)^k/bc;p)}
{\ta (a(w/q)^k/b, a(w/r)^k/c, a(w/s)^k/d;p)}\cr
&\quad\qquad\times
\frac{\ta (et^k, au^k, v^{k-n};p)}
{\ta ((u/w)^kv^{k-n}, e(t/w)^kv^{k-n}/a, e(tu/w)^k;p)}\bigg]\cr
&=\delta_{n,0}
\end{align}
for $n=0,1,\dots,$ where $a^2v^n=bcde$, $w^3=qrstuv$,
and $\delta_{n,m}$ is the Kronecker delta function.

Setting $w=rs$ and $d=a/c$ in \eqref{ftoa}, we have $e=av^n/b$
and obtain (after doing the simultaneous replacements
$q\mapsto r$, $rs\mapsto rst/q$, $u\mapsto rst/q^2$ and $v\mapsto s$) 
the identity (see \cite[Eq.~(11.6.11)]{GR})
\begin{align}\label{dto1}
&\sum^n_{k=0}\frac{\ta (a(rst/q)^k, br^k/q^k, s^{k-n}/q^k, as^nt^k/bq^k;p)}
{\ta (a, b, s^{-n}, as^n/b;p)}\cr
&\quad\times\frac{(a;rst/q^2, p)_k (b;r,p)_k (s^{-n};s,p)_k (as^n/b;t,p)_k}
{(q;q,p)_k (ast/bq;st/q,p)_k (as^nrt/q;rt/q,p)_k
(brs^{1-n}/q; rs/q,p)_k}q^k\cr
&=\delta_{n,0},
\end{align}
where $n = 0, 1, \ldots\,$, which generalizes \cite[Eq.~(3.6.17)]{GR}.
In particular, if we replace $n$, $a$, $b$, and $k$ 
in the $s = t = q$ case of \eqref{dto1} by $n-m$, $ar^mq^m$,
$br^mq^{-m}$, and $j-m$, respectively, we obtain the orthogonality relation
\begin{equation}\label{kd}
\sum^n_{j=m} a_{nj} b_{jm} = \delta_{n,m}
\end{equation}
with
\begin{equation*}
a_{nj} = \frac{(-1)^{n+j} \ta (ar^jq^j, br^jq^{-j};p)
(arq^n, brq^{-n};r,p)_{n-1}}
{(q;q,p)_{n-j} (arq^n, brq^{-n};r,p)_j
(bq^{1-2n}/a;q,p)_{n-j}},
\end{equation*}
\begin{equation*}
b_{jm} =\frac{(ar^mq^m, br^mq^{-m};r,p)_{j-m}}
{(q, aq^{1+2m}/b;q,p)_{j-m}} \Big(-\frac ab q^{1+2m}\Big)^{j-m} 
q^{2\binom{j-m}2}.
\end{equation*}
This shows that the triangular matrix $A = (a_{nj})$ is the inverse of the 
triangular matrix $B = (b_{jm})$, and yields a theta hypergeometric 
analogue of \cite[Eqs.~(3.6.18)--(3.6.20)]{GR}. It should be noted,
on the contrary, that by replacing $n$ and $k$ in \eqref{ftoa} by $n-m$
and $j-m$ one does not obtain a sum of the form \eqref{kd}.

By proceeding as in the derivation of Eq.~(3.6.22) in \cite{GR},
we find that the latter extends to the
bibasic theta hypergeometric summation formula
\begin{equation}
\ta (a/r, b/r;p) \sum^n_{k=0}
\frac{(aq^k, bq^{-k};r,p)_{n-1}\ta (aq^{2k}/b;p)}
{(q;q,p)_k (q;q,p)_{n-k}(aq^k/b;q,p)_{n+1}} (-1)^k 
q^{\binom k2} = \delta_{n,0}
\end{equation}
for $n = 0, 1, \ldots,$ which when $r=q$ reduces to
\begin{equation*}
_{8}V_{7} (a/b;q/b,  aq^{n-1}, q^{-n}, q^{-2n};q,p)=\delta_{n,0}.
\end{equation*}

Special cases of the summation formula \eqref{csn},
combined with the  argument applied
in \cite[Sec.~3.8]{GR}, can be used to extend
equations (3.8.14) and (3.8.15) of \cite{GR}
to the quadratic theta hypergeometric transformation formulas
\begin{align}
&\sum^n_{k=0}\frac{\ta (acq^{3k};p)}{\ta (ac;p)} 
\frac{(a, b, cq/b;q,p)_k (f, a^2c^2q^{2n+1}/f, q^{-2n}; q^2, p)_k}
{(cq^2, acq^2/b, abq;q^2,p)_k
(acq/f, f/acq^{2n}, acq^{2n+1}; q,p)_k} q^k\cr
&=\frac{(acq;q,p)_{2n} (ac^2q^2/bf,abq/f;q^2,p)_n}
{(acq/f;q,p)_{2n} (abq, ac^2q^2/b;q^2,p)_n}\cr
&\quad \times {}_{12}V_{11} (ac^2/b; f, ac/b, c, cq/b, cq^2/b, 
a^2c^2q^{2n+1}/f, q^{-2n};q^2,p)
\end{align}
and
\begin{align}
&\sum^{2n}_{k=0}\frac{\ta (acq^{3k};p)}
{\ta (ac;p)}\frac{(d, f, a^2c^2q/df; q^2,p)_k
(a, cq^{2n+1}, q^{-2n};q,p)_k}
{(acq/d, acq/f, df/ac;q,p)_k
(cq^2, aq^{1-2n}, acq^{2n+2};q^2, p)_k} q^k\cr
&=\frac{(acq, acq/df;q,p)_n (acq^{1-n}/d, acq^{1-n}/f;q^2,p)_n}
{(acq/d, acq/f;q,p)_n (acq^{1-n}, acq^{1-n}/df;q^2,p)_n}\cr
&\quad\times{} _{12}V_{11} (acq^{-2n-1}; c, d, f, a^2c^2q/df, 
aq^{-2n-1}, q^{1-2n}, q^{-2n};q^2,p)
\end{align}
for $n = 0, 1, \ldots\,$; see Thms.~4.2 and 4.7 in Warnaar~\cite{W}.

Also of interest is the special case of \eqref{m0} that is obtained by setting
$w\mapsto rs$, $c\mapsto a/d$, and  $f\mapsto a/d$ (hence $g\to ad/be$),
which after the simultaneous replacements $q\mapsto r$,
$rs\mapsto rst/q$, $u\mapsto rst/q^2$, $v\mapsto s$, $a\mapsto ad$,
and $e\mapsto ad^2/bc$
gives the identity (see also \cite[Eq.~(11.6.9)]{GR})
\begin{align}\label{m00}
&\sum^n_{k=0}
\frac{\ta (ad(rst/q)^k, br^k/dq^k, cs^k/dq^k, adt^k/bcq^k;p)}
{\ta (ad, b/d, c/d, ad/bc;p)}\cr
&\quad\times\frac{(a;rst/q^2,p)_k (b;r,p)_k (c;s,p)_k
(ad^2/bc;t,p)_k}
{(dq;q,p)_k (adst/bq;st/q,p)_k (adrt/cq;rt/q,p)_k
(bcrs/dq;rs/q,p)_k}q^k\cr
&=\frac{\ta(a, b, c, ad^2/bc;p)}{d\,\ta(ad, b/d, c/d, ad/bc;p)}\cr
&\quad\times\frac{(arst/q^2;rst/q^2,p)_n (br;r,p)_n (cs;s,p)_n
(ad^2t/bc;t,p)_n}
{(dq;q,p)_n (adst/bq;st/q,p)_n (adrt/cq;rt/q,p)_n
(bcrs/dq;rs/q,p)_n}\cr
&\quad -\frac{\ta(d, ad/b, ad/c, bc/d;p)}
{d\,\ta(ad, b/d, c/d, ad/bc;p)}.
\end{align}

Just as in the derivation in Gasper~\cite{G} of the quadbasic transformation 
formula in \cite[Ex.$\,$3.21]{GR}, one can extend
indefinite summation formulas (such as in \eqref{sumf} and \eqref{m00})
to transformation formulas by applying the identity
\begin{equation*}
\sum^n_{k=0} \lm_k \sum^{n-k}_{j=0} \Lambda_j = \sum^n_{k=0} \Lambda_k 
\sum^{n-k}_{j=0} \lm_j,
\end{equation*}
which follows by a reversing the order of summation. For example, by taking 
$\lm_k$ to be the $k$th term in the series in \eqref{sumf} and $\Lambda_k$ 
to be this term with $a_k$, $b_k$, $c_k$, $d_k$, and $p$ replaced
by $A_k$, $B_k$, $C_k$, $D_k$, and $P$, respectively, we obtain the 
rather general transformation formula
\begin{align}\label{gtf}
&\sum^n_{k=0}\frac{\ta (a_k, a_k/b_kc_k, a_k/b_kd_k, a_k/c_kd_k;p)}
{\ta (a_k/b_k c_kd_k, a_k/d_k, a_k/c_k, a_k/b_k;p)}\cr
&\quad\times \prod^{k-1}_{j=0}\frac{\ta (b_j, c_j, d_j, a^2_j/b_jc_j d_j;p)}
{\ta (a_j/b_j, a_j/c_j, a_j/d_j, b_jc_jd_j/a_j; p)}\cr
&\quad\times \bigg\{ 1 - \prod^{n-k}_{j=0}
\frac{\ta (B_j, C_j, D_j, A^2_j/B_jC_jD_j;P)}
{\ta (A_j/B_j, A_j/C_j, A_j/D_j, B_jC_jD_j/A_j;P)}\bigg\}\cr
&=\sum^n_{k=0}\frac{\ta (A_k, A_k/B_kC_k, A_k/B_kD_k, A_k/C_k D_k;P)}
{\ta (A_k/B_k C_k D_k, A_k/D_k, A_k/C_k, A_k/B_k;P)}\cr
&\quad\times \prod^{k-1}_{j=0}
\frac{\ta (B_j, C_j, D_j, A^2_j/B_j C_j D_j;P)}
{\ta (A_j/B_j, A_j/C_j, A_j/D_j, B_j C_j D_j/A_j;P)}\cr
&\quad\times \bigg\{ 1 - \prod^{n-k}_{j=0}
\frac{\ta (b_j, c_j, d_j, a^2_j/b_jc_jd_j;p)}
{\ta (a_j/b_j, a_j/c_j, a_j/d_j, b_jc_jd_j/a_j;p)}\bigg\}.
\end{align}

The special case of \eqref{gtf} that is obtained  by using \eqref{m00}
instead of \eqref{sumf} is
\begin{align}\label{ex28}
&\sum^n_{k=0}\frac{\theta (ad (rst/q)^k, br^k/dq^k,
cs^k/dq^k, adt^k/bcq^k;p)}{\theta (ad, b/d, c/d, ad/bc;p)}\cr
&\quad \times\frac{(a; rst/q^2,p)_k (b;r,p)_k (c;s,p)_k (ad^2/bc;t,p)_k}
{(dq;q,p)_k (adst/bq;st/q,p)_k (adrt/cq;rt/q,p)_k (bcrs/dq;rs/q,p)_k} q^k\cr
&\quad\times\bigg(\frac{\ta(A, B, C, AD^2/BC;P)(Q^{-n}/D;Q,P)_k
(B(Q/ST)^n/AD;ST/Q,P)_k}
{D\,\ta(AD, B/D, C/D, AD/BC;P)((Q^2/RST)^n/A;RST/Q^2,P)_k}\cr
&\qquad\times\frac{(C(Q/RT)^n/AD;RT/Q,P)_k (D(Q/RS)^n/BC;RS/Q,P)_k}
{(R^{-n}/B;R,P)_k (S^{-n}/C;S,P)_k (BCT^{-n}/AD^2;T,P)_k} \cr
&\qquad -\frac{\ta(D, AD/B, AD/C, BC/D;P)(DQ;Q,P)_n (ADST/BQ;ST/Q,P)_n }
{D\,\ta(AD, B/D, C/D, AD/BC;P)(ARST/Q^2;RST/Q^2,P)_n (BR;R,P)_n}\cr
&\qquad\times\frac{(ADRT/CQ;RT/Q,P)_n (BCRS/DQ;RS/Q,P)_n}
{(CS;S,P)_n (AD^2T/BC;T,P)_n}\bigg)\cr
&=\frac{(arst/q^2;rst/q^2,p)_n (br;r,p)_n
(cs;s,p)_n (ad^2t/bc;t,p)_n}
{(dq;q,p)_n (adst/bq;st/q,p)_n (adrt/cq;rt/q,p)_n
(bcrs/dq;rs/q,p)_n}\cr
&\quad\times\frac{(DQ;Q,P)_n (ADST/BQ;ST/Q,P)_n}
{(ARST/Q^2;RST/Q^2,P)_n (BR;R,P)_n}\cr
&\quad\times\frac{(ADRT/CQ;RT/Q,P)_n(BCRS/DQ;RS/Q,P)_n}
{(CS;S,P)_n (AD^2T/BC;T,P)_n}\cr
&\times\sum^n_{k=0}\frac{\theta (AD (RST/Q)^k, BR^k/DQ^k,
CS^k/DQ^k, ADT^k/BCQ^k;P)}{\theta (AD, B/D, C/D, AD/BC;P)(DQ;Q,P)_k } Q^k\cr
&\quad \times\frac{(A; RST/Q^2,P)_k (B;R,P)_k(C;S,P)_k (AD^2/BC;T,P)_k}
{(ADST/BQ;ST/Q,P)_k (ADRT/CQ;RT/Q,P)_k (BCRS/DQ;RS/Q,P)_k}\cr
&\quad\times\bigg(\frac{\ta(a, b, c, ad^2/bc;p)
(q^{-n}/d;q,p)_k (b(q/st)^n/ad;st/q,p)_k}
{d\,\ta(ad, b/d, c/d, ad/bc;p)
((q^2/rst)^n/a;rst/q^2,p)_k (r^{-n}/b;r,p)_k }\cr
&\qquad\times\frac{(c(q/rt)^n/ad;rt/q,p)_k (d(q/rs)^n/bc;rs/q,p)_k}
{(s^{-n}/c;s,p)_k (bct^{-n}/ad^2;t,p)_k} \cr
&\qquad -\frac{\ta(d, ad/b, ad/c, bc/d;p)(dq;q,p)_n (adst/bq;st/q,p)_n}
{d\,\ta(ad, b/d, c/d, ad/bc;p)(arst/q^2;rst/q^2,p)_n (br;r,p)_n}\cr
&\qquad\times\frac{(adrt/cq;rt/q,p)_n (bcrs/dq;rs/q,p)_n}
{(cs;s,p)_n (ad^2t/bc;t,p)_n}\bigg).
\end{align}

The $d,D\to 1$ special case of \eqref{ex28} is
\begin{align}
&\sum^n_{k=0}\frac{\theta (a (rst/q)^k, br^kq^{-k},
cs^kq^{-k}, at^k/bcq^k;p)}{\theta (a, b, c, a/bc;p)}\cr
&\quad \times\frac{(a; rst/q^2,p)_k (b;r,p)_k (c;s,p)_k (a/bc;t,p)_k}
{(q;q,p)_k (ast/bq;st/q,p)_k (art/cq;rt/q,p)_k (bcrs/q;rs/q,p)_k}\cr
&\quad\times\frac{(Q^{-n}; Q,P)_k (B(Q/ST)^n/A; ST/Q,P)_k
(C(Q/RT)^n/A; RT/Q,P)_k}
{((Q^2/RST)^n/A; RST/Q^2,P)_k (R^{-n}/B;R,P)_k (S^{-n}/C;S,P)_k}\cr
&\quad \times\frac{((Q/RS)^n/BC;RS/Q,P)_k}{(BC/AT^n;T,P)_k} q^k\cr
&=\frac{(arst/q^2;rst/q^2,p)_n (br;r,p)_n (cs;s,p)_n (at/bc;t,p)_n}
{(q;q,p)_n (ast/bq;st/q,p)_n (art/cq;rt/q,p)_n (bcrs/q;rs/q,p)_n}\cr
&\quad \times\frac {(Q;Q,P)_n (AST/BQ; ST/Q,P)_n}
{(ARST/Q^2;RST/Q^2,P)_n (BR;R,P)_n}\cr
&\quad \times\frac {(ART/CQ; RT/Q,P)_n (BCRS/Q; RS/Q,P)_n}
{(CS; S,P)_n (AT/BC;T,P)_n}\cr
&\quad \times \sum^n_{k=0}\frac{\theta (A(RST/Q)^k, BR^kQ^{-k},
CS^kQ^{-k}, AT^k/BCQ^k;P)}{\theta (A, B, C, A/BC;P)}\cr
&\qquad\times\frac{(A;RST/Q^2,P)_k (B;R,P)_k}
{(Q;Q,P)_k (AST/BQ; ST/Q,P)_k}\cr
&\qquad \times\frac{(C;S,P)_k (A/BC;T,P)_k}
{(ART/CQ;RT/Q,P)_k (BCRS/Q;RS/Q,P)_k}\cr
&\qquad \times\frac{(q^{-n};q,p)_k (b(q/st)^n/a; st/q,p)_k}
{((q^2/rst)^n/a; rst/q^2,p)_k}\cr
&\qquad\times\frac{(c(q/rt)^n/a;rt/q,p)_k ((q/rs)^n/bc;rs/q,p)_k}
{(r^{-n}/b;r,p)_k (s^{-n}/c;s,p)_k (bc/at^n; t,p)_k}Q^k,
\end{align}
for $n = 0, 1, \ldots$ . For $s=t=q$ and $S=T=Q$ this reduces to
the elliptic quadbasic transformation formula
\begin{align}
&\sum^n_{k=0}\frac{\theta (a r^k q^k, br^k q^{-k};p)}
{\theta (a, b;p)}\frac{(a, b; r, p)_k (c, a/bc;q,p)_k}
{(q, aq/b;q,p)_k (ar/c, bcr;r,p)_k} \cr
&\quad \times\frac{(CR^{-n}/A, R^{-n}/BC; R,P)_k (Q^{-n}, BQ^{-n}/A;Q,P)_k}
{(Q^{-n}/C, BCQ^{-n}/A; Q, P)_k (R^{-n}/A, R^{-n}/B;R, P)_k}q^k\cr
&=\frac{(ar, br;r,p)_n (cq, aq/bc;q,p)_n (Q, AQ/B;Q,P)_n (AR/C, BCR; R, P)_n}
{(q, aq/b;q,p)_n (arc, bc/r;r,p)_n (AR, BR; R,P)_n (CQ, AQ/BC;Q,P)_n}\cr
&\quad \times \sum^n_{k=0}\frac{\theta (AR^kQ^k, BR^kQ^{-k};P)}
{\theta (A, B; P)}\frac{(A, B; R, P)_k (C, A/BC; Q, P)_k}
{(Q, AQ/B; Q, P)_k (AR/C, BCR; R, P)_k} \cr
&\quad\times\frac{(cr^{-n}/a, r^{-n}/bc;r, p)_k (q^{-n}, bq^{-n}/a;q,p)_k}
{(q^{-n}/c, bcq^{-n}/a; q,p)_k (r^{-n}/a, r^{-n}/b;r,p)_k}Q^k,
\end{align}
which is an extension of the second  identity in \cite[Ex.$\,$3.21]{GR}
(see also \cite[Ex.$\,$11.25]{GR}).
If we now set $R=Q=r=q$, we obtain
the following transformation formula for a
``split-poised'' theta hypergeometric series
\begin{align}\label{splitpoised1}
&\sum^n_{k=0}\frac{\theta (aq^{2k};p)}{\theta (a;p)}
\frac{(a, b, c, a/bc;q,p)_k}{(q, aq/b, aq/c, bcq;q,p)_k}\cr
&\quad \times\frac{(q^{-n}, B/Aq^n, C/Aq^n, 1/BCq^n;q,p)_k}
{(1/Aq^n, 1/Bq^n, 1/Cq^n, BC/Aq^n;q,p)_k} q^k\cr
&=\frac{(aq, bq, cq, aq/bc, Aq/B, Aq/C, BCq; q,p)_n}
{(Aq, Bq, Cq, Aq/BC, aq/b, aq/c, bcq;q,p)_n}\cr
&\quad \times \sum^n_{k=0}
\frac{\theta (Aq^{2k};p)}{\theta (A;p)}
\frac{(A, B, C, A/BC;q,p)_k}{(q, Aq/B, Aq/C, BCq;q,p)_k}\cr
&\qquad\times\frac{(q^{-n}, b/aq^n, c/aq^n, 1/bcq^n;q,p)_k}
{(1/aq^n, 1/bq^n, 1/cq^n, bc/aq^n;q,p)_k} q^k
\end{align}
for $n=0,1,\ldots,$ which is an extension of the
transformation formula for  a split-poised $_{10}\phi_9$ series
given in \cite[Ex.$\,$3.21]{GR}.
This formula may also be written as a transformation formula for a
split-poised $_{12}E_{11}$ series, see \eqref{splitpoised}.

We will now use \eqref{dto1} to derive multibasic extensions of the
Fields and Wimp,  Verma, and Gasper expansion formulas in
\cite[Eqs.~(3.7.1)--(3.7.3) \& (3.7.6)--(3.7.9)]{GR}, and multibasic
theta hypergeometric extensions of \cite[Eqs.~(3.7.6)--(3.7.8)]{GR}.
Let $a = \gamma (rst/q)^j$ and $b = \s (r/q)^j$ in \eqref{dto1}
and replace the summation index $k$ by $n-k$.
For $j,n = 0, 1, \ldots\,$, we assume that
$B_n(p)$ and $C_{j,n}$ are complex numbers such that $C_{j,0} = 1$ 
and the sequence $\big\{B_n(p)\big\}$
has finite support when $p \ne 0$.
Then, for $j = 0, 1, \ldots,$  
\begin{align}\label{knj}
B_j(p) x^j&=\sum^\ty_{n=0}\frac{\ta(\gamma\s^{-1}(st)^{n+j};p)
(\gamma\s^{-1} (st)^j;st/q,p)_n}
{\ta(\gamma\s^{-1} (st)^j;p) (s;s,p)_n}\cr
&\ \ \times \frac{(\gamma rs^jtq^{-1};rt/q,p)_j
(\s rs^{1-j}q^{-1};rs/q,p)_j}
{(\gamma rs^{n+j}tq^{-1};rt/q,p)_j(\s rs^{1-n-j}q^{-1};rs/q,p)_j}\cr
&\ \ \times s^{\binom{n+1}2}q^{-\binom{n+1}2-nj}
B_{j+n}(p) C_{j,n} x^{j+n}\delta_{n,0} \cr
&=\sum^\ty_{k=0}\sum^\ty_{n=j}\frac{\ta (\gamma (rst/q)^n, \s (r/q)^n, 
\gamma\s^{-1}(st)^{n+k}, \gamma\s^{-1} s^{n+k}t^nq^{j-n}, s^{-k}q^{j-n};p)}
{\ta (s^{j-n-k};p) (s;s,p)_k (q;q,p)_n} \cr
&\ \ \times \frac{(\gamma\s^{-1} (st)^{n+1} q^{j-n-1};st/q,p)_{k-1} 
(\gamma \s^{-1}s^{n+k} t^{j+1};t,p)_{n-j-1}}
{(\gamma rs^{n+k}tq^{-1};rt/q,p)_n
(\s rs^{1-n-k}q^{-1};rs/q,p)_n} \cr
&\ \ \times (\gamma rs^jtq^{-1};rt/q,p)_j
(\gamma (rst)^{j+1}q^{-j-2}; rst/q^2,p)_{n-j-1}\cr
&\ \ \times (\s rs^{1-j}q^{-1};rs/q,p)_j
(\s r^{j+1}q^{-j};r, p)_{n-j-1} (q^{-n};q,p)_j \cr
&\ \ \times (-1)^n B_{n+k}(p) C_{j, n+k-j} x^{n+k}
s^{\binom {k+1}2}q^{-\binom {k+1}2+\binom n2+n(1+j-n-k)}
\end{align}
by interchanging sums and setting $n\mapsto n+k-j$
(this extension of \cite[Eq.~(3.7.5)]{GR} corrects \cite[Eq.~(11.6.20)]{GR}).

By multiplying both sides of \eqref{knj} by  $A_j w^j/(q;q,p)_j$
and summing from $j = 0$ to $\ty$ we get that the following multibasic
expansion formula (this corrects \cite[Eq.~(11.6.21)]{GR})
\begin{align}\label{gs1}
&\sum^\ty_{n=0} A_nB_n(p)\frac{(xw)^n}{(q;q,p)_n}\cr
&= \sum^\ty_{n=0}\frac{\ta (\gamma (rst/q)^n, \s (r/q)^n; p)}
{(q;q,p)_n} (-x)^n q^{n+\binom n2}\cr
&\quad\times \sum^\ty_{k=0}\frac{\ta (\gamma\s^{-1}(st)^{n+k};p)}
{(\gamma rs^{n+k}tq^{-1};rt/q,p)_n
(\s rs^{1-n-k}q^{-1};rs/q,p)_n} \frac{B_{n+k}(p) x^k}
{(s;s,p)_k}
s^{\binom {k+1}2}q^{-\binom {k+1}2}\cr
&\quad\times \sum^n_{j=0}\frac{\ta (\gamma\s^{-1} s^{n+k}t^nq^{j-n},
s^{-k}q^{j-n};p)(q^{-n};q,p)_j}
{\ta (s^{j-n-k};p) (q;q,p)_j}\cr
&\quad\times
(\gamma\s^{-1} (st)^{n+1} q^{j-n-1};st/q,p)_{k-1}
(\gamma \s^{-1}s^{n+k} t^{j+1};t,p)_{n-j-1} \cr
&\quad\times (\gamma rs^jtq^{-1};rt/q,p)_j
(\gamma (rst)^{j+1}q^{-j-2}; rst/q^2,p)_{n-j-1}\cr
&\quad\times (\s rs^{1-j}q^{-1};rs/q,p)_j
(\s r^{j+1}q^{-j};r, p)_{n-j-1} A_j C_{j, n+k-j} w^j q^{n(j-n-k)},
\end{align}
which reduces to \cite[Eq.~(3.7.6)]{GR} by letting $p=0$ and then
setting $r = p$ and $s = t = q$.

If we set $r = s = t = q$ and $C_{j,m}\equiv 1$ in \eqref{gs1}
we obtain an expansion formula that is equivalent to the
following extension of \cite[Eq.~(3.7.7)]{GR}
(which corrects a slight misprint in \cite[Eq.~(11.6.22)]{GR})
\begin{align}\label{gs2}
&\sum^\ty_{n=0} A_nB_n(p)\frac{(xw)^n}{(q;q,p)_n} =
\sum^\ty_{n=0}\frac{(\s, \gamma q^{n+1}/\s, \al, \be; q,p)_n}
{(q, \gamma q^n;q,p)_n} \Big(\frac x{\s}\Big)^n\cr
&\quad\times \sum^\ty_{k=0}\frac{\ta (\gamma q^{2n+2k}/\s;p)
(\gamma q^{2n}/\s, \s^{-1}, \al q^n, \be q^n;q,p)_k}
{\ta (\gamma q^{2n}/\s;p) (q, \gamma q^{2n+1};q,p)_k} B_{n+k}(p) x^k\cr
&\quad\times \sum^n_{j=0}\frac{(q^{-n}, \gamma q^n;q,p)_j}
{(q, \gamma q^{n+1}/\s, q^{1-n}/\s, \al, \be; q,p)_j}
A_j (wq)^j,
\end{align}
where, as previously, it is assumed that $\big\{B_n(p)\big\}$
has finite support when $p\ne 0$.
Clearly, one cannot let $\s \to \ty$ in \eqref{gs2} to obtain an
extension of \cite[Eq.~(3.7.3)]{GR} that holds when $p\ne 0$.

Corresponding to the $q$-analogue of the Fields and Wimp expansion
formula displayed in \cite[Eq.~(3.7.8)]{GR},  \eqref{gs2}
gives the rather general theta hypergeometric
expansion formula (see \cite[Eq.~(11.6.23)]{GR})
\begin{align}\label{gs3}
&\sum^\ty_{n=0}\frac{(a_R, c_T;q,p)_n}
{(q, b_S, d_U; q,p)_n} A_nB_n(p) (xw)^n\cr
&=\sum^\ty_{n=0}\frac{(c_T, e_K, \s, \gamma q^{n+1}/\s;q,p)_n}
{(q, d_U, f_M, \gamma q^n;q,p)_n} \Big(\frac x{\s}\Big)^n\cr
&\quad\times \sum^\ty_{k=0}\frac{\ta (\gamma q^{2n+2k}/\s;p)
(\gamma q^{2n}/\s, \s^{-1}, c_T q^n, e_K q^n;q,p)_k}
{\ta (\gamma q^{2n}/\s;p)
(q, \gamma q^{2n+1}, d_U q^n, f_M q^n;q,p)_k} B_{n+k}(p) x^k\cr
&\quad\times \sum^n_{j=0}\frac{(q^{-n}, \gamma q^n, a_R, f_M;q,p)_j}
{(q, \gamma q^{n+1}/\s, q^{1-n}/\s, b_S, e_K; q,p)_j} A_j (wq)^j,
\end{align}
where we used the contracted notation that was used in
\cite[Eq.~(3.7.8)]{GR}, and, in order to avoid convergence problems, it is
assumed that $\big\{B_n\big\}$ has finite support when $p\ne 0$.

By using \eqref{ftoa} instead of its special case 
\eqref{dto1} and proceeding as
above, one can derive even more general extensions of the
multibasic Fields and Wimp, Verma, and Gasper expansions.
Since they are rather lengthy and do not seem to be of any particular
interest at this time, we will not give them here.

\end{document}